\newcommand{\ncm}{\newcommand}
\ncm{\aut}{auto\-mor\-phi\-sm}
\ncm{\cstar}{$C^{*}$-algebra}
\ncm{\cstars}{$C^{*}$-algebras}
\ncm{\ra}{\mbox{$\rightarrow$}}
\ncm{\al}{\mbox{$\alpha $}}
\ncm{\Ad}{\mbox{\rm Ad}}
\ncm{\ep}{\mbox{$\epsilon > 0$}}
\ncm{\ol}{\overline}
\ncm{\MR}{M. R\o{}rdam}
\ncm{\vo}{}
\ncm{\ch}{}
\ncm{\CMP}{Comm. Math. Phys.}
\ncm{\add}{}
\ncm{\BE}{\begin{eqnarray*}}
\ncm{\EE}{\end{eqnarray*}}
\ncm{\id}{\mbox{\rm id}}
\ncm{\iso}{\cong}

\newtheorem{theo}{Theorem}[section]
\newtheorem{cor}[theo]{Corollary}
\newtheorem{lem}[theo]{Lemma}
\newtheorem{prop}[theo]{Proposition}
\newtheorem{remark}[theo]{Remark}
\newtheorem{definition}[theo]{Definition}

\newenvironment{rem}{\begin{remark} \rm}{\end{remark}}
\newenvironment{pf}{{\it Proof.}}{\vspace{3mm}}
\newenvironment{df}{\begin{definition} \rm}{\end{definition}}

\documentstyle[12pt]{article}
\input{amssym.def}
\input{amssym.tex}

\newfont{\mine}{msbm10 at 12pt}

\def\Z{{\hbox{\mine Z}}}
\def\R{{\hbox{\mine R}}}
\def\N{{\hbox{\mine N}}}
\def\C{{\hbox{\mine C}}}
\def\Q{{\hbox{\mine Q}}}

\title{Trace scaling automorphisms of certain stable AF~algebras}
\author{ \normalsize David E. Evans and
         Akitaka Kishimoto
  }
\date{}
\begin{document}
\maketitle

\begin{abstract}
Trace scaling \aut s of a stable AF algebra with dimension group totally
ordered are outer conjugate if the scaling factors are the same
(not equal to one).
\end{abstract}

\section{Introduction}
This is a continuation of \cite{EEK}, where we showed a UHF version
of a well-known result of A. Connes \cite{C2} that trace scaling
\aut s of the AFD type II$_{\infty}$ factor with the same non-trivial scale
are outer conjugate with each other. In this paper we show the same
result for stable AF algebras with totally ordered dimension group.

The key idea remains the same as in \cite{EEK} and hence as in \cite{C2}:
Define and prove a Rohlin property for such automorphisms and analyse them
using this property. We are now familiar with the unital case
(see \cite{BKRS,BEK,KK,K1,K2,K3}). What we did in \cite{EEK} is to
evade non-unital \cstar s and deal with {\em partial} unital endomorphisms
of unital (UHF) algebras.
What we do here is to define a suitable Rohlin property 
for \aut s of non-unital \cstar s and prove it. 
We will define it by borrowing an idea due to R\o{}rdam
\cite{R}, where corner endomorphisms are treated, and prove it 
by using an argument in \cite{K1,K2}, 
where \aut s of unital AF algebras are
treated. Our main contribution is to find a passage from the {\em non-unital}
case to the {\em unital} case in proving the Rohlin property, 
which is done in Section 2. 

In Section 3 we shall show that our definition of Rohlin property is
the {\em right} one, i.e., this is at least strong enough to prove
the stability or 1 cohomology property \cite{C2,HO1,HO2}. Note that
it is this property that we actually need.

Another idea used in \cite{C2,EEK}, a technique envolving tensor
products, is no longer applicable here. For example if $A$ is an AF algebra
with totally ordered dimension group and $A\otimes A$ has the same
property, then $A$ must be UHF. However
by using an intertwining argument we shall show that any two
\aut s \al\ and $\beta$ with the Rohlin property of an AF algebra are
outer conjugate if $\al_*$ and $\beta_*$ are equal as automorphisms of
the dimension group. More precisely, the conclusion is that
for any \ep\ there is an \aut\ $\sigma$ of $A$ and a unitary $U$ in
$A+\C 1$ (or in $A$ if $A$ is unital) 
such that $\|U-1\|<\epsilon$ and
$\al=\Ad\,U\circ\sigma\circ\beta\circ\sigma^{-1}$.
In general we cannot take $1$ for $U$ or cannot conclude conjugacy
of $\alpha$ and $\beta$.
This extends Theorem 2 of \cite{K2} and Theorem 7 of \cite{EEK}.
(Note that even in the UHF case 
the conclusion is stronger than the one in \cite{EEK}.)

Let $\lambda\neq1$ be a positive number and let $G_{\lambda}$ 
be the subgroup of
\R\ generated by $\lambda^n,\, n\in\Z$. If $A$ is the stable AF algebra
whose dimension group is $G_{\lambda}$ and \al\ is an automorphism of
$A$ such that $\al_*$ acts on the dimension group by
multiplication by $\lambda$, then from the above result the crossed
product $A\rtimes_{\alpha}\Z$ depends only on $\lambda$. 
But we have now a more general theorem in this direction:
it follows from \cite{R,Kir,P} that $A\rtimes_{\alpha}\Z$, being
a purely infinite simple \cstar, is isomorphic to a stable Cuntz algebra.
 
\section{Rohlin property}
\setcounter{theo}{0}
Let $A$ be a non-unital \cstar\ and let \al\ be an \aut\ of $A$. 
We assume that $A$ has an approximate unit consisting of projections.
Based on \cite{R} we define a Rohlin property for \al\ as
follows:

\ncm{\F}{\mbox{${\cal F}$}}
\begin{df}\label{2.1} 
The \aut\ \al\ has the Rohlin property if for any $k\in\N$ there are
positive integers $k_1,\ldots,k_m\geq k$ satisfying the following
condition: For any projections $E,e$ in $A$ , any unitary $U$ in
$A+\C1$, any finite subset \F\ of $A_E=EAE$, and \ep\ with
$$
 e\leq E,\ \Ad\,U\circ\al(e)\leq E,\ e\in\F,\ \Ad\,U\circ\al(e)\in\F,
$$
there exists a family $\{e_{i,j};\ i=1,\ldots,m,\ j=0,\ldots,k_i-1\}$
of projections in $A$ such that
\BE
&& \sum_i\sum_j e_{i,j}=E,\\
&& \|\Ad\,U\circ\al(e_{i,j}e)-e_{i,j+1}\Ad\,U\circ\al(e)\|<\epsilon,\\
&& \|[x,e_{i,j}]\|<\epsilon,
\EE
for $i=1,\ldots,m,\ j=0,\ldots,k_i-1$ and $x\in\F$ where 
$e_{i,k_i}=e_{i,0}$. The projections $\{e_{i,j}\}$ will be called a set
of Rohlin towers.
\end{df}

If we apply the same definition to a unital \cstar\ $A$, then
the Rohlin property for the unital case \cite{K1,K2}, where
$E=1=e$ is preassumed, implies the present definition.
(We just have to cut down by $E$ a set of Rohlin towers obtained for
$E=1=e$ which almost commutes with $E,\, e$ and use functional calculus
to get the desired set of Rohlin towers.)
The following is an easy consequence whose proof we omit:

\begin{prop}\label{A}
Suppose that there is an increasing sequence $\{P_n\}$ of projections
in $A$ such that $\|P_nx-x\|\ra0$ for any $x\in A$ and $\al(P_n)=P_n$.
Then \al\ has the Rohlin property if and only if the restriction of
\al\ to $A_{P_n}=P_nAP_n$ has the Rohlin property for any $n$.
\end{prop}

We are, however, interested in the situation where the above proposition
does not apply.

Let $A$ be a simple stable AF algebra and let \al\ be an automorphism
of $A$. Let $\{A_n\}$ be an increasing sequence of finite-dimensional
subalgebras of $A$ such that the union $\cup_n A_n$ is dense in $A$.

\begin{rem}\label{REMARK}
In this situation for any \ep\ there is a subsequence $\{n_k\}$ of
positive integers and a unitary $U\in A+1$ such that
$\|U-1\|<\epsilon$ and $\Ad\,U\circ\al(A_k)\subset A_{k+1},\ 
(\Ad\,U\circ\al)^{-1}(A_k)\subset A_{k+1}$ for any $k$.
This can be proved by using the following fact inductively:
If $B$ is a finite-dimensional subalgebra of $A$ and \ep,
there is an $n\in\N$ and a unitary $U\in A+1$ such that
$\|U-1\|<\epsilon$ and $\Ad\,U(B)\subset A_n$.
\end{rem}

Hence by slightly perturbing \al\ and passing to a subsequence of
$\{A_n\}$ we may assume that
$\al^{-1}(A_n)\subset A_{n+1},\ \al(A_n)\subset A_{n+1}$
for any $n$. 
We fix a nonzero projection $E\in A_1$.

\ncm{\alp}{\hbox{$\Ad\,U\circ\al$}}
Let $e$ be a projection in $\cup_k A_k$ and $U$ a unitary in
$\cup_k A_k + \C 1$ such that 
$e\leq E$ and $\Ad\, U\circ\al(e)\leq E$.
Since $A_E=EAE$ is simple there is a $k\in\N$ such that
$e,\ \Ad\,U\circ\al(e),\ U\, \in A_k+\C1$ and the multiplication by
$e$ (resp. $f=\alp(e)$) induces an isomorphism of
$(A\cap A_k')_E$ onto $(A\cap A_k')_e$ 
(resp. $(A\cap A_k')_f$) or the central support of $e$ (resp. $f$) in
$(A_k)_E$ is $E$.
We define a homomorphism $\phi(\al,e,U)$ of
$(A\cap A_{k+1}')_E$ into $(A\cap A_k')_E$ by
$$
\phi(\al,e,U)(x)f=\alp(xe)=\al(x)f.
$$
This is indeed well-defined: Since
$$
[\alp(xe),b]=\alp([xe,(\alp)^{-1}(b)])=0
$$
for $b\in (A_k)_f=fA_kf$, we have that
$$
\alp(xe)\in (A\cap A_k')_f
=A_f \cap ((A_k)_f)'.
$$
\ncm{\phietc}{\hbox{$\phi(\al,e,U)$}}
We note that \phietc\ is essentially independent of $e$ and $U$ in the sense
that if $\phi(\al,e_1,U_1)$ is another one then 
$\phietc=\phi(\al,e_1,U_1)$ on the common domain $D$.
Because if $e_1\leq e$ and $U_1=U$, this follows since
$$
(\phietc(x)-\phi(\al,e_1,U_1)(x))\alp(e_1)=0
$$ 
for any $x\in D$; if $e=e_1$, this follows since
$$
(\phietc(x)-\Ad(UU_1^*)\circ\phi(\al,e,U_1)(x))\alp(e)=0
$$
for any $x\in D$ and $UU_1^{*}$ commutes with $\phi(\al,e,U_1)(x)$;
and if $e_1=\Ad\,V(e)$ and $U_1=U\al(V^*)$ with $V$ a unitary in
$\cup_n A_n +\C1$, this follows since 
$\alp(e)=\Ad\,U_1\circ\al(e_1)$.

Thus we denote by $\tilde{\al}$ the homomorphism induced by these
\phietc; this is a homomorphism of $(A\cap A_{k+1}')_E$ into
$(A\cap A_k')_E$ for some $k$. The same computation shows that
$\tilde{\al}$ maps $(A\cap A_{n+1}')_E$ into $(A\cap A_n')_E$
for $n\geq k$. Since $\phi(\al^{-1},\alp(e),\al^{-1}(U^*))$ is
well-defined if so is \phietc, $\tilde{\al^{-1}}$ is defined
at least on $(A\cap A_{k+2}')_E$, and satisfies that
$\tilde{\al}\circ\tilde{\al^{-1}}=\id,\ \tilde{\al^{-1}}\circ\tilde{\al}=\id$
on $(A\cap A_{k+2}')_E$. A similar computation shows that
$(\tilde{\al})^n=\tilde{\al^n}$ on $(A\cap A_l')_E$ for a sufficiently
large $l$ (depending on $n$).

Let $\omega$ be a free ultrafilter on \N\ and let $A_E^{\omega}$
be the quotient \cstar\ of $l^{\infty}(\N,A_E)$ by the ideal
$I_{\omega}=\{ (x_n)\, |\, \lim_{n\ra\omega}\|x_n\|=0\}$.
Embedding $A_E$ into $l^{\infty}(\N, A_E)$ and so into $A_E^{\omega}$
as constant functions, let $A_{E\omega}=A_E^{\omega}\cap A_E'$.

Let $x=(x_n) \in A_{E\omega}$. Then we can find an increasing
sequence $\{k_n\}$ in \N\ and $x_n'\in (A\cap A_{k_n}')_E$
such that $k_n\ra \infty$ and $\lim_{n\ra\omega}\|x_n-x_n'\|=0$.
We define a homomorphism $\tilde{\al}_{\omega}$ of
$A_{E\omega}$ into itself by
$$
  \tilde{\al}_{\omega}(x)=(\tilde{\al}(x_n')).
$$
This is indeed easily checked to be well-defined. In the same way
we can define $\tilde{\beta}_{\omega}$ for $\beta=\al^{-1}$ and show that
$\tilde{\al}_{\omega}\circ\tilde{\beta}_{\omega}=\id$
and $\tilde{\beta}_{\omega}\circ\tilde{\al}_{\omega}=\id$.
Thus $\tilde{\al}_{\omega}$ is an automorphism of $A_{E\omega}$.

Let $\tau$ be a densely-defined non-zero lower semi-continuous trace
on $A$ and assume that $\tau$ is unique up to a constant multiple.
Since $\tau\circ\al$ is again such a trace, there is a
$\lambda>0$ such that $\tau\circ\al=\lambda\tau$. 
We normalize $\tau$ by $\tau(E)=1$. We can define
a state $\tau_{\omega}$ on $A_{E\omega}$ by
$\tau_{\omega}(x)=\lim \tau(x_n)$ for $x=(x_n)\in A_{E\omega}$.
Note that $\tau_{\omega}$ is tracial. We shall show that
$\tau_{\omega}$ is invariant under $\tilde{\al}_{\omega}$.

Let $x=(x_n)\in A_{E\omega}$, where $x_n\in (A\cap A_{k_n})_E'$
for some non-decreasing sequence $\{k_n\}$ with $k_n\ra \infty$.
If $\tilde{\al}=\phi(\al,e,U)$ on $(A\cap A_{k+1}')_E$, 
$\{P_i\,;\ i=1,\ldots,N\}$ is the set of minimal projections
in the center of $EA_{k+1}E$, and $k_n>k+1$, then
\BE
 \tau(\tilde{\al}(x_n))&=& \sum_{i=1}^N \tau(P_i\tilde{\al}(x_n))\\
   &=& \sum\frac{\tau(P_i)}{\tau(P_i\Ad\,U\circ\al(e))}
       \tau(P_i\Ad\,U\circ\al(x_ne))\\
   &=& \sum\frac{\tau(P_i)}{\tau(P_i\al(e))}\tau\circ\al(\al^{-1}(P_i)ex_n)\\
   &=& \sum\frac{\tau(P_i)}{\tau(\al^{-1}(P_i)e)}\tau(\al^{-1}(P_i)ex_n)
\EE
which, when $n$ is large, is almost equal to
$$
  \sum\frac{\tau(P_i)}{\tau(\al^{-1}(P_i)e)}\tau(\al^{-1}(P_i)e)\tau(x_n)=\tau(x_n)
$$
since $\tau$ is factorial. Thus we obtain that 
$\tau_{\omega}\circ\tilde{\al}_{\omega}=\tau_{\omega}$.

Without loss of generality we assume, from now on, that $\tilde{\al}$ is 
defined as $\phi(\al,e,1)$ on $(A\cap A_2')_E$, i.e., $e,\, \al(e)\in A_1,\
e\leq E,\, \al(e)\leq E$, and the central supports of $e$ and $\al(e)$
in $(A_1)_E$ are $E$.
\ncm{\cR}{\mbox{${\cal R}$}}

The above argument carries over to the weak closure 
$\cR$ of $\pi_{\tau}(A_E)$. Note that $\cR^{\omega}$ is defined as the 
quotient of $l^{\infty}(\N,\cR)$ by
$$ 
I= \{(x_n)\, |\, \lim_{n\ra\omega}\|x_n\|_{\tau}=0\}
$$
where $\|a\|_{\tau}=\tau(a^*a)^{1/2}$, and $\tau$ is regarded as the tracial 
state on \cR\ induced from $\tau$ on $A_E$. Since 
$\tau\circ\tilde{\al}|(A\cap A_2')_E$ is equivalent to
$\tau|(A\cap A_2')_E$, $\tilde{\al}$ extends to a homomorphism of
$\cR\cap \pi_{\tau}(A_{2E})'$ into 
$\cR\cap \pi_{\tau}(A_{1E})'$, and for 
$x\in l^{\infty}(\N,\cR\cap \pi_{\tau}(A_{2E})')$,
$\lim_{n\ra\omega}\|x_n\|_{\tau}=0$ if and only if 
$\lim\|\tilde{\al}(x_n)\|_{\tau}=0$.
In this way we have the automorphism $\tilde{\al}_{\omega}$ of
$\cR_{\omega}=\cR^{\omega}\cap\cR'$ induced by $\tilde{\al}$ which satisfies 
that $\tau_{\omega}\circ\tilde{\al}_{\omega}=\tau_{\omega}$.

\begin{lem}\label{2.2}
Suppose that $\tau\circ\al=\lambda\tau$ with $\lambda\neq1$.
If $\cR=\pi_{\tau}(A_E)^{-w},\ \cR_{\omega},\ \tilde{\al}_{\omega},\ 
\tau_{\omega}$ etc. are as above, then any non-zero power of 
$\tilde{\al}_{\omega}$ is properly outer.
\end{lem}
\begin{pf}
The proof is similar to the proofs of Lemmas 1 and 2 of \cite{EEK}.
Denote by $\ol{\al}$ the \aut\
of $\pi_{\tau}(A)''$ induced by \al.
Fix $n\geq2$ and let $B=\pi_{\tau}(A_{nE})\subset \cR$.
For a non-zero projection $f\in \cR\cap B'$ we assert that
$$
 \inf\{\|p\tilde{\al}(p)\|\,;\ 0\neq p=p^*=p^2\in \cR\cap B',\ 
                  p\leq f\}=0.
$$

Suppose that the above infimum is positive, say $\delta>0$.
We may suppose that $f$ is in one factor direct summand of 
$\cR\cap B'$. Let $f_1$ be a minimal projection in 
$B$ such that $f_1f\neq0$. Then any projection $\tilde{p}\in\cR$
with $\tilde{p}\leq ff_1$ is of the form
$\tilde{p}=f_1p$ with $p$ a projection in 
$f(\cR\cap B')f$. Hence for any $z\in\cR$ we have
$$
\inf\{\|pf_1z\ol{\al}(f_1)\ol{\al}(p)\|\ ;\ 0\neq p=p^*=p^2 
          \in\cR\cap B',\ p\leq f \}=0,
$$
since $\ol{\al}$ is an outer \aut\ of $\pi_{\tau}(A)''$.
There is a finite set $\{V_1,\ldots,V_k\}$ of unitaries in $B$ such that
$ \sum_iV_if_1V_i^*$ is the central support of $f_1$ in $B$.
Since
$$
\sum_{i,j}V_ipf_1V_i^*\ol{\al}(V_j)\ol{\al}(f_1)\ol{\al}(p)\ol{\al}(V_j^*)
=p\ol{\al}(p)
$$
and
$$
\|p\tilde{\al}(p)\|=\|p\tilde{\al}(p)\ol{\al}(e)\|=\|p\ol{\al}(p)\ol{\al}(e)\|
=\|p\ol{\al}(p)\|,
$$
there exist $i,j$ such that
$$
  \|pf_1V_i^*\ol{\al}(V_j)\ol{\al}(f_1)\ol{\al}(p)\|\geq \delta/k^2,
$$
which is a contradiction.
By using this we can show that $\al_{\omega}$ is properly outer
(cf. Lemma 2 of \cite{EEK} and \cite{C2}).  

By applying the same argument to $\tilde{\al}^n$, we obtain that
$(\tilde{\al}_{\omega})^n$ is properly outer for any $n\neq0$.
\end{pf}

Let $D=A_E$ and suppose that $D$ has a unique tracial state.
Let $B$ be a finite-dimensional subalgebra of $D$. We say 
$x\in D$ is {\em independent} of $B$ if $x\in B'$ and
$$
\tau(xy)=\tau(x)\tau(y),\ y\in B.
$$

\begin{lem}\label{2.3}
Let $\{P_1,\ldots,P_N\}$ be the set of minimal central projections of $B$.
Then $x\in D\cap B'$ is independent of $B$ if and only if
$$
  \tau(xP_i)=\tau(x)\tau(P_i),\ i=1,\ldots,N.
$$
\end{lem}
\begin{pf}
It suffices to show the {\em if} part. Since
$P_iDP_i\cong P_iBP_i\otimes (D\cap B')P_i$, it follows that
for $y\in P_iBP_i$
$$
  \frac{\tau(xy)}{\tau(P_i)}=\frac{\tau(xP_i)}{\tau(P_i)}\frac{\tau(y)}{\tau(P_i)}.
$$
Hence $\tau(xy)=\tau(x)\tau(y)$. This completes the proof.
\end{pf}

\begin{lem}\label{2.4}
Suppose that the dimension group $K_0(A)$ is totally ordered 
and identified with a subgroup of \R\ and
that $K_0(A)=\lambda K_0(A)$ for some $\lambda\neq 1$.
For any central sequence $\{f_k\}$ of projections in $A_E$ there is a
central sequence $\{f_k'\}$ of projections in $A_E$ such that
$f_k'\leq f_k$, $\tau(f_k-f_k')\ra 0$, and for any $n$ there is a 
$k_n$ satisfying that $f_k'$ is independent of $A_{nE}$ for
any $k\geq k_n$.
\end{lem}
\begin{pf}
We may suppose that there is a $k_n$ such that $f_k\in (A\cap A_n')_E$
for $k\geq k_n$. Let $\{P_i^{(n)}\}$ be the set of minimal central
projections in $A_{nE}$. We can find an $l_n\geq k_n$ such that
$$
|\frac{\tau(P_i^{(n)}f_k)}{\tau(P_i^{(n)})}-\tau(f_k)|<\frac{1}{n}
$$
for $k\geq l_n$ and all $i$. We will then define $f_k',\ k\in\{l_n,\ldots,
l_{n+1}-1\}$ as follows: Let $g_k\in \cup_{n\in \Z}\lambda^n\Z$
be such that $\max\{0,\gamma-1/n\} <g_k\leq\gamma$ where
$$
\gamma=\min_i\frac{\tau(P_i^{(n)}f_k)}{\tau(P_i^{(n)})}.
$$
Choose a subprojection $q_{k,i}$ of $P_i^{(n)}f_k$ in
$(A\cap A_n')_{P_i^{(n)}f_k}$ such that
$\tau(q_{k,i})/\tau(P_i^{(n)})=g_k$, and let
$f_k'=\sum_i q_{k,i}$. This is possible because, when $(A_n)_{P_i^{(n)}}$
is isomorphic to the $m_i\times m_i$ matrix algebra, 
$$
 K_0((A\cap A_n')_{P_i^{(n)}})=m_iK_0(A_{P_i^{(n)}})
$$
and when $[P_i^{(n)}]=1$,
$K_0((A\cap A_n')_{P_i^{(n)}})$ contains $\cup_{n\in\Z}\lambda^n\Z$
for any $i$.
\end{pf}

\begin{lem}\label{2.5}
Suppose that $K_0(A)$ is totally ordered. Let $p\in (A\cap A_n')_E$
be a projection independent of $(A_n)_E$. Then if $n-m\geq 1$, 
$p,\tilde{\al}(p),\ldots,\tilde{\al}^m(p)$ are all equivalent
in $(A\cap A_{n-m}')_E$.
\end{lem}
\begin{pf}
Let $\{P_i\}$ be the set of minimal central projections in 
$EA_{n-1}E$. Then
\BE
 \tau(\tilde{\al}(p)P_i)&=& \frac{\tau(P_i)}{\tau(P_i\al(e))}
                              \tau(\tilde{\al}(p)P_i\al(e))\\
  &=&\frac{\tau(P_i)}{\tau(\al^{-1}(P_i)e)}\tau(\al^{-1}(P_i)ep)\\
  &=&\tau(P_i)\tau(p)=\tau(pP_i),
\EE
since $\al^{-1}(P_i)e\in A_n$.
Hence $p$ is equivalent to $\tilde{\al}(p)$ in
$(A\cap A_{n-1}')_E$. We just repeat this procedure.
\end{pf}

\begin{theo}\label{2.6}
Let $A$ be a stable AF algebra such that
$K_0(A)$ is totally ordered and let \al\ be an \aut\ of $A$ such that
$\tau\circ\al=\lambda\tau$ where $\tau$ is a trace on $A$ (unique up to
constant multiple) and $\lambda\neq 1$. Then \al\ has the Rohlin property.
\end{theo}
\begin{pf}
Let $E,e$ be projections in $A$ and $U$ a unitary in $A+\C1$ such that
$e\leq E$ and $\Ad\, U\circ\al(e)\leq E$. Let $\{A_n\}$ be an increasing 
sequence of finite-dimensional subalgebras of $A$ such that the union
$\cup_n A_n$ is dense in $A$ and $E,\, e,\, \Ad\,U\circ\al(e)\in A_1$.
By taking $\Ad\,U\circ\al$ instead of \al\ we now assume that $U=1$.
For any $\delta>0$ we find a unitary $V\in A+1$ such that
$\|V-1\|<\delta$ and, by passing to a subsequence of $\{A_n\}$,
$\Ad\,V\circ\al(A_n)\subset A_{n+1},\ (\Ad\,V\circ\al)^{-1}(A_n)\subset
A_{n+1}$ for any $n$
(Remark \ref{REMARK}). 
By perturbing $V$ if necessary we may further assume that
$\Ad\,V\circ\al(e)=\al(e)$. By taking a sufficiently small $\delta>0$
we may take $\Ad\,V\circ\al$ for \al. Thus we have the following situation:
There exists an increasing sequence $\{A_n\}$ of finite-dimensional
subalgebras of $A$ such that the union $\cup_n A_n$ is dense in $A$,
$\al(A_{n})\subset A_{n+1},\ \al^{-1}(A_n) \subset A_{n+1},\ 
E,e,\al(e)\in A_1$,
and $e, \al(e)\leq E$. The problem is to find a set of Rohlin towers
as specified in Definition \ref{2.1}. But this follows from
Lemmas \ref{2.2}, \ref{2.4}, and \ref{2.5} 
based on the arguments given in \cite{K1,K2}
because we just have to prove the (ordinary) Rohlin property
for $\tilde{\al}$. Here is an outline. Since $\cR_{\omega}$ is
a finite von Neumann algebra \cite[2.2.1]{C2}, we know
by Lemma 2.3 and \cite[1.2.5]{C2} that $\tilde{\al}_{\omega}$
on $\cR_{\omega}$ satisfies a Rohlin property. Then by approximating a
Rohlin tower in $\cR_{\omega}$ by projections in $A_E$ we show that
the {\em partial} endomorphism $\tilde{\al}$ of $A_E$ satisfies an
approximate Rohlin property (Lemmas 2.5 and 2.6 and \cite{K2}).
But this suffices to conclude the Rohlin property 
by \cite[Proof of 2.1]{K1}.
\end{pf}

\begin{rem}\label{2.7}
In the situation of the previous theorem let \al\ be an \aut\ of $A$
such that $\al_*=\id$ and any non-zero power is not weakly inner in
the tracial representation. Then \al\ has the Rohlin property.
(See Proposition \ref{A} or Theorem 4.1 of \cite{K2}.)
\end{rem}

\begin{rem}\label{2.8}
In the above theorem we can make the Rohlin property more specific:
In Definition \ref{2.1} we may take $\{k,k+1\}$ for
$\{k_1,\ldots,k_m\}$. This follows because of Lemma \ref{2.5} (see \cite{K2}).
\end{rem}

\section{Stability}
\setcounter{theo}{0}
\begin{theo}\label{3.1}
Let $A$ be a (non-unital) AF algebra and let \al\ be an \aut\ of $A$
with the Rohlin property.
Let $\epsilon>0$ and let $B_1$ be a finite-dimensional subalgebra of $A$.
Then there is a finite-dimensional subalgebras $B_2$ of $A$ such that
for any unitary $U\in A\cap B_2'+1$ there is a unitary 
$V\in A\cap B_1'+1$ with
$\|U-V\al(V^*)\|<\epsilon$.
\end{theo}

The following argument works if $A$ is unital or non-unital; if $A$ is
unital, we should regard $1$ as the unit of $A$, otherwise $1$ as
the unit adjoined to $A$.

Suppose that there is an increasing sequence $\{A_n\}$ 
of finite-dimensional subalgebras of $A$ such that
$\cup_n A_n$ is dense in $A$ and
$\al^{-1}(A_{n})\subset A_{n+1},\ \al(A_n)\subset A_{n+1}$ for any $n$.
By Remark \ref{REMARK}
the above theorem is an easy consequence of:

\begin{lem}
Under the above assumption let $U$ be a unitary in $A\cap A_n'+1$.
Then for each $k\in\N$ there is a unitary $V$ in
$A\cap A_{n-2k}'+1$ such that 
$\|U-V\al(V^*)\|<4/k$, where $A_m=\{0\}$ for a non-positive $m$.
\end{lem}

\begin{pf}  
Let $U$ be a unitary in $A\cap A_n'+1$. We may further suppose that
there is an $m>n$ such that $U\in A_m+1$. Let $F$ be the identity of
$A_m$ and let
\BE
&&  E=\sup\{\al^m(F);\ -1\leq m\leq 2k+1\},\\
&&  e=\sup\{\al^m(F);\ -1\leq m\leq 2k\},
\EE
which are projections in $A_{m+2k+1}$ with
that $e, \al(e)\leq E$. For
$U_j=U\al(U)\cdots\al^{j-1}(U)$ with $j\geq0$, we have that
if $0\leq j\leq k$,
$$ U_j(1-\al(e))=1-\al(e).
$$
For $E,e$ and $\{k,k+1\}$ we find a set of Rohlin towers
$\{e_{1,0},\ldots,e_{1,k-1},\,e_{2,0},\ldots,e_{2,k}\}$
in $(A\cap A_{m+2k+1}')_E$. Let $W_t^{(1)},\, W_t^{(2)}$ be paths of unitaries
in $A_{m+k}\cap A_{n-k}'+1$ such that
\BE
&& W_0^{(i)}=1,\\
&& W_1^{(1)}=U_k,\ W_1^{(2)}=U_{k+1},\\
&& \|W_s^{(i)}-W_t^{(i)}\|\leq\pi|s-t|,\ \ s,t\in [0,1],\\
&& W_t^{(i)}(1-e_0)=1-e_0,\ \ \ t\in [0,1],
\EE
where $e_0=\sup\{\al^m(F);\ 0\leq m\leq k\}$.
Set
$$
V=\sum_{j=0}^{k-1}U_j\al^j(W_{1-j/(k-1)}^{(1)})e_{1,j}+
  \sum_{j=0}^k U_j\al^j(W_{1-j/k}^{(2)})e_{2,j} +1-E,
$$
which is a unitary in $A\cap A_{n-2k}'+1$.
Since $\al^j(W_t^{(i)})(1-\al(e))=1-\al(e)$ for $0\leq j\leq k$, we have
$$ V(1-\al(e))=1-\al(e),\ \al(V)(1-\al(e))=1-\al(e).
$$
Hence it follows that 
\BE
 V\al(V^*)\al(e)&
\approx&\sum_{j=0}^{k-1} U_{j+1}\al^{j+1}(W_{1-(j+1)/(k-1)}^{(1)})
             \al^{j+1}(W_{1-j/(k-1)}^{(1)})^*\al(U_j^*)e_{1,j+1}\al(e)\\
&&+\sum_{j=0}^k U_{j+1}\al^{j+1}(W_{1-(j+1)/k}^{(2)})\al^{j+1}(W_{1-j/k}^{(2)})^*\al(U_j^*)e_{2,j+1}\al(e),
\EE
where the $k-1$'th summand in the first summation should be understood as
$$
 U_0W_1^{(1)}\al^k(W_0^{(1)})^*\al(U_{k-1}^*)e_{1,0}\al(e)=Ue_{1,0}\al(e)
$$
and the $k$'th term in the second as
$$
U_0W_1^{(2)}\al^{k+1}(W_0^{(2)})^*\al(U_k^*)e_{2,0}\al(e)=Ue_{2,0}\al(e).
$$
Hence we have that
$$
\|V\al(V^*)-U\|\leq \pi/k + \epsilon(2k+1),
$$
where $\epsilon>0$ is a small number 
depending on the Rohlin towers we used. This completes
the proof.
\end{pf}

\section{Outer conjugacy}
\setcounter{theo}{0}
\begin{theo}\label{4.1}
Let $A$ be an AF algebra and let \al\ and $\beta$ be \aut s of $A$ with the 
Rohlin property. If $\al_*=\beta_*$ on $K_0(A)$, then for any \ep\
there is an \aut\ $\sigma$ of $A$ and a unitary $U$ in $A+1$ such that
$\al=\Ad\,U\circ\sigma\circ\beta\circ\sigma^{-1},\
\|U-1\|<\epsilon$, and $\sigma_*=\id$.
\end{theo}
\begin{pf}
Note that $A$ can be either unital or non-unital.

Let \ep\ and let $\{x_n\}$ be a dense sequence in the unit ball of $A$.
We shall construct inductively sequences $\{A_n\},\ \{B_n\}$ of
finite-dimensional subalgebras of $A$, sequences $\{u_n\},\ \{v_n\}$
of unitaries in $A+1$ such that
\begin{enumerate}
\item $A_n \ni_{1/n} x_1,\ldots,x_n,\ A_n\supset A_{n-1}$,
\item $B_n\supset A_n,\ B_n\ni v_n$,
\item $\beta_{2n}|A_{2n+1}=\Ad\,u_{2n+1}\circ\al_{2n-1}|A_{2n+1}$,
\item $\al_{2n-1}|A_{2n}=\Ad\,u_{2n}\circ\beta_{2n-2}|A_{2n}$,
\item $\|u_{2n+1}-v_{2n+1}\al_{2n-1}(v_{2n+1}^*)\|<2^{-2n-1}\epsilon,\
         v_{2n+1}\in B_{2n-1}'$,
\item $\|u_{2n}-v_{2n}\beta_{2n-2}(v_{2n}^*)\|<2^{-2n}\epsilon,\
         v_{2n}\in B_{2n-2}'$,
\item for any unitary $U\in A\cap \al_{2n-1}(A_{2n})'+1$ there exists a
    a unitary $V\in A\cap B_{2n-1}'+1$ such that $\|U-V\al_{2n-1}(V^*)\|<
    2^{-2n-1}\epsilon$,
\item for any unitary $U\in A\cap \beta_{2n}(A_{2n+1})'+1$ there exists a
   unitary $V\in A\cap B_{2n}'+1$ such that
   $\|U-V\beta_{2n}(V^*)\|<2^{-2n-2}\epsilon$,
\end{enumerate}
where $A_n\ni_{\delta}x$ means that there is a $y\in A_n$ with
$\|x-y\|<\delta$ and
\BE
&& A_0=\{0\},\\
&& \al_{2n+1}=\Ad\,u_{2n+1}\circ\al_{2n-1},\ \ \al_{-1}=\al,\\
&& \beta_{2n}=\Ad\,u_{2n}\circ\beta_{2n-2},\ \ \beta_0=\beta.
\EE

We first construct $A_1$ according to (1).
Having constructed 
$$
A_1,\ldots,A_{2n+1},\ B_1,\ldots,B_{2n},\
u_1,\ldots,u_{2n},\ v_1,\ldots,v_{2n},
$$ 
we proceed as follows:
We choose $u_{2n+1}$ according to (3). Since $\al_{2n-1}|A_{2n}=
\beta_{2n}|A_{2n}$ from (4) and the definition of $\beta_{2n}$ above, 
it follows that 
$u_{2n+1}\in \al_{2n-1}(A_{2n})'$ and so by (7) that there is a unitary
$v_{2n+1}\in A\cap B_{2n-1}'+1$ satisfying (5). We may assume that there
is a $B_{2n+1}$ satisfying (2) (for $2n+1$ in place of $n$). 
Having defined $B_{2n+1}$ we define
$A_{2n+2}$ satisfying (1) and (7) (by using Theorem \ref{3.1}) and choose $u_{2n+2}$ according to
(4). Since $u_{2n+2}\in \beta_{2n}(A_{2n+1})'$, we define
$v_{2n+2}$ satisfying (6) by using (8) and assume that there is a
$B_{2n+2}$ satisfying (2). We define $A_{2n+3}$ satisfying (1) and (8).
This completes the induction.

We note that the union $\cup_n A_n$ is dense in $A$ and
define \aut s $\sigma_n$ of $A$ by
\BE
&& \sigma_{2n}=\Ad(v_{2n}v_{2n-2}\cdots v_2),\\
&& \sigma_{2n+1}=\Ad(v_{2n+1}v_{2n-1}\cdots v_1),
\EE
and define
\BE 
&& \tilde{\sigma}_0=\lim_n \sigma_{2n},\\
&& \tilde{\sigma}_1=\lim_n \sigma_{2n+1}.
\EE
Since $v_{n-2},v_{n-4},\ldots \in B_{n-2}$,
$v_n\in B_{n-2}'$, and $\cup_n B_n$ is dense, they are well-defined. We let
$$
w_{2n+1}= u_{2n+1}\al_{2n-1}(v_{2n+1})v_{2n+1}^*,\ \ 
w_{2n}= u_{2n}\beta_{2n-2}(v_{2n})v_{2n}^* 
$$
and define unitaries
$w_n'\in A+1$ by
\BE
&& w_{2n}'=w_{2n}\Ad\, v_{2n}(w_{2n-2})\Ad(v_{2n}v_{2n-2})(w_{2n-4})\cdots
        \Ad(v_{2n}\cdots v_4)(w_2),\\
&& w_{2n+1}'=w_{2n+1}\Ad\,v_{2n+1}(w_{2n-1})\cdots \Ad(v_{2n+1}\cdots v_3)(w_1).
\EE
Since $\|w_n-1\|<2^{-n}\epsilon$, both
$\{w_{2n}'\}$ and $\{w_{2n+1}'\}$ converge, say to $\tilde{w}_0$ and
$\tilde{w}_1$ respectively.
Then $\tilde{w}_i$'s  are unitaries in $A+1$ such that
$\|\tilde{w}_i-1\|<\epsilon$.

Since $\al_{2n-1}|A_{2n}=\beta_{2n}|A_{2n}$, we have that
$$
\Ad\,w_{2n-1}'\circ\sigma_{2n-1}\circ\al\circ\sigma_{2n-1}^{-1}|A_{2n}=
  \Ad\,w_{2n}'\circ\sigma_{2n}\circ\beta\circ\sigma_{2n}^{-1}|A_{2n},
$$
which implies that
$$
\Ad\,\tilde{w}_1\circ\tilde{\sigma}_1\circ\al\circ\tilde{\sigma}_1^{-1}=
\Ad\,\tilde{w}_0\circ\tilde{\sigma}_0\circ\beta\circ\tilde{\sigma}_0^{-1}.
$$
This completes the proof.
\end{pf}

\begin{cor}
Let $A$ be a stable AF algebra such that $K_0(A)$ is totally ordered
and let $\al, \beta$ be \aut s of $A$ such that $\tau\circ\al=\lambda\tau$
and $\tau\circ\beta=\lambda\tau$ where $\tau$ is a trace on $A$ and
$\lambda\neq 1$. Then for any \ep\ there are an \aut\ $\sigma$ of $A$
and a unitary $U$ in $A+1$ such that 
$\|U-1\|<\epsilon,\ \sigma_*=\id$, and 
$\al=\Ad\,U\circ\sigma\circ\beta\circ\sigma^{-1}$.
\end{cor}
\begin{pf}
This follows from Theorems \ref{2.6} and \ref{4.1}.
\end{pf}

\begin{rem}
In the above corollary the exact conjugacy $\al=\sigma\circ\beta
\circ\sigma^{-1}$
for some \aut\ $\sigma$ of $A$ cannot be expected in general.
For example if $A$ is a prime AF algebra such that
$A\iso A\otimes{\cal K}$, where ${\cal K}$ is the compact operators on
$l^2(\Z)$, and \al\ is an \aut\ of $A$ such that 
$\al^n$ is properly outer for any $n\neq0$, let $\al_1$ be the \aut\ of $A$
defined as $\al\otimes\gamma$ through $A\otimes{\cal K}\iso A$
where $\gamma=\Ad\,U$ and $U$ is the shift unitary on $l^2(\Z)$.
Then $\al_1$ satisfies that 
for any $x,y \in A$,
$$
\|\al_1^n(x)y\|\ra 0
$$
as $n\ra\infty$. 
This property is preserved by conjugacy but not by outer conjugacy.
(By \cite{K0} there exists a faithful $\al_1$-covariant irreducible
representation of $A$; by using Kadison's transitivity theorem in
this irreducible representation it follows that for any \ep\ there
are an $x\in A$ and a unitary $U\in A+1$ such that $0\leq x \leq 1$,
$\|x\|=1$, $\|U-1\|<\epsilon$, and 
$\|(\Ad\,U\circ\al_1)^n(x)x\|=1$.)
\end{rem}

\begin{rem}
Let $\lambda\neq1$ be a positive number and let
$$
G_{\lambda}=\cup_{n\in \Z}\Z\lambda^n.
$$
If $\{1,\lambda,\lambda^2,\ldots,\}$ are linearly independent over $\Q$
then the quotient $G_{\lambda}/(1-\lambda)G_{\lambda}$
is isomorphic to $\Z$ and otherwise
if $\{f\in \Z[t]\ |\ f(\lambda)=0\}=p(t)\Z[t]$ with some
$p(t)\in \Z[t]$, then $G_{\lambda}/(1-\lambda)G_{\lambda}\iso\Z/p(1)\Z$.
If $A$ is the stable AF algebra with dimension group $G_{\lambda}$
and \al\ is an \aut\ of $A$ with $\al_*=\lambda$, then the crossed product
$A\rtimes_{\alpha}\Z$ has $G_{\lambda}/(1-\lambda)G_{\lambda}$ as
$K_0$ and $\{0\}$ as $K_1$ by the Pimsner-Voiculescu
exact sequence \cite{Bl}. Hence $A\rtimes_{\alpha}\Z$ is isomorphic to a stable
Cuntz algebra ${\cal O}_n\otimes {\cal K}$ where $n$ is either
finite or infinite \cite{R,Kir,P,Cu1,Cu2}.
\end{rem}

\medskip
\small
\begin{flushright}
\begin{tabular}{l}
Department of Mathematics\\
University of Wales\\
Swansea SA2 8PP, UK\\
\hfill\\
Department of Mathematics\\
Hokkaido University\\
Sapporo 060 Japan
\end{tabular}
\end{flushright}
\end{document}